\documentclass[a4paper,10pt]{amsart}
\usepackage{amssymb}

\usepackage[T1]{fontenc}
\usepackage[utf8]{inputenc}
\usepackage[german,english]{babel}

\usepackage[babel,german=guillemets]{csquotes} %F�r Literaturverwaltung
\usepackage[backend=bibtex8,style=authoryear]{biblatex} %F�r Literaturverwaltung
\bibliography{emptheorien2.bib}%F�r Literaturverwaltung

\usepackage{paralist}
\usepackage[bottom]{footmisc}
 
\begin{document}
\pagestyle{plain}

\setlength{\skip\footins}{3pt}
\renewcommand{\footnoterule}{\rule[3pt]{2cm}{0.3pt}}
\setlength{\footnotesep}{5mm}
 
\author{Hans Joachim Burscheid $\cdot$ Horst Struve}

\address{Universität zu Köln\\Seminar für Mathematik und ihre Didaktik\\Gronewaldstr. 2, 50931 Köln\\Tel.: 0221--470--4750, Fax.: 0221--470--4985\\e--mail: h.burscheid@uni--koeln.de, h.struve@uni--koeln.de}

\title{Empirische Theorien im Kontext der Mathematikdidaktik}

\selectlanguage{german}
\begin{abstract}
An zwei methodologisch unterschiedlich bearbeiteten Beispielen wird gezeigt, wie der Begriff der empirischen Theorie (in strukturalistischer Darstellung) genutzt werden kann, mathematikdidaktisch relevanten Inhalten eine präzise Form zu geben. 
\end{abstract}

\maketitle

\vspace{-1cm}

\selectlanguage{english}
\begin{abstract}
At two examples dealt with in methodologically different ways it will be pointed out how the concept of an empirical theory (in the sense of the Structuralists) can be useful to specify contents relevant to maths didactics.
\end{abstract}

\thispagestyle{empty}

\section{\textbf{Einleitung}}
Wie nicht zuletzt das vor wenigen Jahren erschienene Themenheft des ZDM "Comparing, Combining, Coordinating -- Networking Strategies for Connecting Theoretical Approaches" [Prediger u. a. 2008] belegt, besteht in der Mathematikdidaktik ein zunehmendes Interesse an Theoriebildung. Man könnte etwas bescheidener formulieren: Ein zunehmendes Interesse, die gewonnenen Erkenntnisse in einer verbindlichen Form zu organisieren. In diesem Sinne ist der folgende Beitrag zu verstehen.

An zwei methodologisch unterschiedlich bearbeiteten Beispielen wollen wir zeigen, welche Möglichkeiten der Begriff der empirischen Theorie --- und ihrer formalen Darstellung --- bieten, um mathematikdidaktische Inhalte in eine strukturierte Form zu bringen. Bevor wir jedoch Inhalte thematisieren, erscheint es uns zweckmäßig, den Leser\footnote{und natürlich auch die Leserin} vorab davon zu überzeugen, daß es auch in der Mathematikdidaktik hilfreich sein kann, etablierte Formen der Darstellung zu übernehmen und auch dort formale Darstellungen zu wählen, wo sie bisher eher unüblich sind. Daß dies möglich und sinnvoll ist, hat als Ursache, daß sich etliche Inhalte der Mathematikdidaktik als sog. empirische Theorien darstellen lassen. Den Begriff \emph{empirische Theorie} kann man als Oberbegriff aller Theorien verstehen, die ein Stück WELT beschreiben, wie es z.B. die Theorien der Naturwissenschaftler tun. Eine etablierte \emph{Form zur Darstellung} empirischer Theorien ist die strukturalistische Metatheorie, die von Wolfgang Stegmüller und seinem Kreis entwickelt wurde [Stegmüller 1973, 1986], [Balzer 1982]. Sie ist nicht nur in den Naturwissenschaften sondern auch in Wirtschaftswissenschaft, Psychologie, Genetik und anderen wissenschaftlichen Disziplinen anwendbar [Balzer u.a. 2000]. Ihre Formalisierung erfolgt durch ein mengentheoretisches Prädikat, wie man es von den Axiomatisierungen der Mathematik kennt, z.B. (G; $\circ$) heißt eine \emph{Gruppe}, wenn ... . Ein wesentliches Merkmal dieser Darstellungsform ist, daß sie den Aufbau der Theorie wiedergibt. Die Darstellung gliedert sich in drei Stufen. Auf der ersten Stufe werden als sog. \emph{partielle Modelle} der Theorie die empirischen Gegebenheiten formuliert, über die die Theorie Aussagen machen kann. Zu den partiellen Modellen gehören insbesondere die sog. \emph{intendierten Anwendungen} der Theorie, die Gegebenheiten, zu deren Beschreibung und Erklärung die Theorie entwickelt wird. Damit gehen charakteristische Anwendungen in die Darstellung der Theorie ein. So kann diese nicht von ihrem Kontext getrennt werden. Die nächste Stufe --- die der \emph{potentiellen Modelle} --- dient der Ergänzung der Sprache, in der die Theorie formuliert wird. Während auf der ersten Stufe noch alle Begriffe reale Referenzen haben oder aus bekannten (früher erworbenen) Theorien stammen, bedarf die neu zu formulierende Theorie neuer Begriffe, will sie neues Wissen vermitteln. Es sind die bzgl. der in den partiellen Modellen formulierten empirischen Gegebenheiten \emph{theoretischen} Begriffe, solche, die dort keine Referenzen haben, die erst durch die neu zu formulierende Theorie eine Bedeutung erhalten. Auf der dritten Stufe werden schließlich die Axiome formuliert, die die \emph{Modelle} der Theorie ausmachen, die Strukturen, die die intendierten Anwendungen der Theorie erklären. Da eine empirische Theorie keine universelle Anwendung kennt sondern nur einzelne Anwendungen, wird abschließend durch die sog. \emph{Querverbindung} festgelegt, welche der potentiellen Modelle der Theorie von dieser tatsächlich erfaßt werden. \bigskip 

Das erste Beispiel, das wir im folgenden aufgreifen, geht auf die Diskussion der Zahlaspekte zurück und verwendet nur den begrifflichen Aufbau der Darstellung einer empirischen Theorie, wie er gerade geschildert wurde. Ähnlich wie die Aspekte einer natürlichen Zahl in ihrer Gesamtheit das Verfügen über den Begriff der natürlichen Zahl ausmachen, zeigte Robert W. Lawler mit seinen Mikrowelten, wie sich das erste Zahl -- und Operationsverständnis aus dem Zusammenwirken unterschiedlicher, sich zunächst gegenseitig ausschließender Erfahrungsbereiche (in der Terminologie von Bauersfeld) entwickelt [Lawler 1981]. Er unterscheidet Instrumentale Welten, Serien -- Welten und Konforme Welten. Im Gegensatz zu den Instrumentalen Welten, die empirisch verifizierbar sind, haben die nur skizzierten Serien -- Welten und Konformen Welten den Charakter mentaler Modelle. Unseres Wissens sind sie nie mit Beispielen belegt worden, die nicht aus dem Lawlerschen Kontext stammten. Indem wir die Begriffe der Serien -- Welt und der Konformen Welt durch Relationen zwischen Instrumentalen Welten ersetzen, gelingt es, die Lawlersche Idee so zu rekonstruieren, daß an die Stelle von Vagheiten in der Formulierung von Serien -- Welt und Konformer Welt präzise Interpretationen treten, und das Konzept auch auf andere Beispiele übertragbar wird. \medskip

\noindent \emph{Bemerkung}: Wenn wir von einer Rekonstruktion der Lawlerschen Idee sprechen, so meint Rekonstruktion hier nicht eine Direktinterpretation mit der Absicht, die Lawlerschen Gedanken nachzuzeichnen und möglichst genau seine Intentionen zu treffen. Gemeint ist vielmehr eine sog. \emph{rationale Rekonstruktion}. Der Begriff geht auf Stegmüller zurück. Dieser bezeichnet damit eine Analyse historischer Beiträge zur Philosophie [Stegmüller 1967]. Eine rationale Rekonstruktion ist durch drei Prinzipien bestimmt:
\begin{quote}
\begin{small}
"(1) Die Theorie muß in solcher Form dargeboten werden, daß die Darstellung mit den Grundideen des betreffenden Philosophen in Einklang bleibt.
(2) Soweit wie möglich soll die Theorie mittels \emph{präziser} Begriffe dargestellt werden.\\
(3) Sie soll als \emph{konsistente Theorie} entwickelt werden, falls dies möglich ist (d.h. falls sich nicht alle rationalen Deutungen, welche die Forderungen (1) und (2) erfüllen, als inkonsistent erweisen)"
[ebd., S. 2].
\end{small}
\end{quote}
Speziell zu (3) führt Stegmüller aus: \\
\begin{small}
"Das Prinzip (3) wird in all jenen Fällen relevant, in denen verschiedene, einander widersprechende rationale Interpretationen möglich sind. In Bezug auf eine Rekonstruktion ergibt es keinen Sinn zu fragen, ob sie eine wahre Deutung liefere, da sie auf jeden Fall irgendwie vom Originaltext abweichen muß. Sie kann nur \emph{mehr oder weniger adäquat sein}" [ebd., S. 4].
\end{small} \smallskip 

Hans Poser hat sich mit der Methode der rationalen Rekonstruktion kritisch auseinandergesetzt. Er beschreibt sie wie folgt [Poser 1971]: \\
\begin{small}
"Das Charakteristikum der Methode ist also ihre dem Ausgang nach systematische Fragestellung, die an historisches Material herangetragen wird, nicht so sehr, um ein Gedankengebäude neu zu durchdringen, als vielmehr, um es --- mit den Mitteln dieser Methode --- für die Behandlung systematischer Probleme dienstbar zu machen" [ebd., S. 69]. \\
\end{small}
 Das Ergebnis einer Untersuchung, die er mit der Methode der rationalen Rekonstruktion durchführte, formuliert er wie folgt: \\
\begin{small}
"So ist die Gewinnung eines systematischen, wenngleich mit einem historischen 'Aufhänger' versehenen Resultates nicht weiter verwunderlich, ja, das legitime Ziel des Ansatzes" [ebd., S. 72]. \\
\end{small}  
Bezüglich des Prinzips (1) sieht er \\
\begin{small}
"die Gefahr, daß nicht ein System und das Verständnis eines Systems, sondern das Selbstverständnis des Übersetzers von einer Problemstellung rational rekonstruiert wird" [ebd., S. 78].
\end{small} \smallskip

Daß es sich im vorliegenden Fall um einen psychologischen und nicht um einen philosophischen Ansatz handelt, der uns erklärungsbedürftig erscheint, dürfte vergleichsweise unerheblich sein. Auf die Frage der Adäquatheit der Rekonstruktion kommen wir noch zurück. \bigskip

Die Intention des zweiten Beispiels ist eine völlig andere. Es dürfte vor allem das Interesse derjenigen finden, die sich um die Identifizierung von Defiziten bemühen, die das Wissen von Schülern aufweist. Zu vielen stoffdidaktischen Themen gibt es inzwischen eine beeindruckende Fülle empirischer Ergebnisse. Diese sind aber häufig nicht erzielt worden, orientiert an einem vorgegebenen theoretischen Rahmen, sondern beziehen sich vielfach auf weitere Untersuchungen des gleichen Problemfeldes ("lokal organisiert"). Es bleibt offen, wie auf der Basis einer solchen Vielfalt eine Theorie entwickelt werden kann, die sich als Grundlage einer Unterrichtskonzeption nutzen läßt. In diesem Punkt glauben wir, kann der Begriff der empirischen Theorie hilfreich sein. Denn die Theorien, die Kinder im Schulunterricht erwerben, sind empirische Theorien, die sich explizit darstellen lassen. Die Kleinschrittigkeit einer formalen Darstellung legt die einzelnen Punkte offen, die für die Theorie konstitutiv sind. Die formale Darstellung erlaubt keine Lücken im Aufbau. Vorliegende empirische Ergebnisse können zur Rechtfertigung einzelner Festlegungen herangezogen und so in die Theorie einbezogen werden. Dies ist ein Vorzug empirischer Theorien. Stoffdidaktische Ausarbeitungen sind mathematischer Natur, und in mathematische Theorieteile lassen sich empirische Ergebnisse formal nicht einbinden. Die formale Darstellung empirischer Theorien ermöglicht es durch ihre Stufung, in einem systematischen --- nicht notwendig genetischen --- Sinne den Wissenserwerb nachzuzeichnen. Liest man die empirische Theorie deskriptiv, so beschreibt sie den Aufbau des Wissens beim Schüler. Liest man sie dagegen präskriptiv / normativ, so kann sie als Vorlage für eine Unterrichtskonzeption dienen, die die in die Theorie aufgenommenen Ergebnisse berücksichtigt. Da die formale Darstellung in geschlossener Form vorliegt, kann mit rationalen Argumenten entschieden werden, was die Konzeption von ihr übernimmt und was nicht. Wir führen dies hier \emph{beispielhaft} am Thema Flächeninhalt durch, d.h. wir erheben nicht den Anspruch, daß die berücksichtigten empirischen Untersuchungen die einzig relevanten sind --- obwohl sie auch nicht zufällig gewählt wurden. Wählt man andere, so kann dies bedeuten, daß die Theorie anders formuliert werden muß, d.h. die zugrundegelegten empirische Ergebnisse \emph{rechtfertigen} die formale Darstellung. Gleichzeitig macht eine solche Darstellung deutlich, in welchen Punkten empirische Ergebnisse zu ihrer Rechtfertigung fehlen. 

Dies wäre ein erster Schritt, eine gewisse Systematik in die Vielfalt der empirischen Ergebnisse zu bringen und diese in sie umfassende Ordnungen einzubinden.

\section{\textbf{Zahlaspekte und Mikrowelten}}
Kommen wir zum ersten Beispiel. Schon unter den ersten Mathematikdidaktikern, die sich nach dem 2. Weltkrieg dem Grundschulunterricht widmeten, waren etliche, die eine mathematische Ausbildung erhalten hatten, die von den Ideen der Bourbakisten der 30er Jahre des vorigen Jahrhunderts beeinflußt war. Zentrale Elemente waren Mengen und Strukturen. Es ist daher nicht weiter erstaunlich, daß die Auffächerung des Zahlbegriffs der Mathematik --- Kardinalzahlen vs. Ordinalzahlen --- in die Grundschuldidaktik übernommen und dort noch weiter aufgeschlüsselt wurde. Außer Kardinal- und Ordinalzahlen unterschied man Maßzahlen, Operatoren, Rechenzahlen, um nur einige Beispiele zu nennen. Generell sprach man von verschiedenen Aspekten des Zahlbegriffs --- hier zunächst nur des Begriffs der natürlichen Zahl. Auf wen die Bezeichnung \emph{Zahlaspekt} zurückgeht, konnten wir nicht ermitteln. Hans Freudenthal verwendet sie schon in seiner Monographie \emph{"Mathematics as an educational task"}, deren Vorwort aus dem Jahre 1970 datiert. Auf eine präzise Definition des Aspektbegriffs wurde allerdings verzichtet. Er wurde nur umgangssprachlich gefaßt, z.B. wie folgt: Ein \emph{Aspekt} der natürlichen Zahlen ist gegeben durch eine Klasse von Anwendungen, in denen die natürlichen Zahlen in gleicher Weise verwendet werden.  \smallskip   
 
In der kognitiven Psychologie begann sich Anfang der 70er Jahre des 20sten Jahrhunderts die Einsicht durchzusetzen, daß  Wissenskonstruktionen bereichspezifisch sind: \emph{frames} (Marvin Minsky), \emph{scripts} (Roger C. Schank), \emph{microworlds} (R. W. Lawler). Heinrich Bauersfeld  machte diese Ideen unter den Mathematikdidaktikern der Bundesrepublik bekannt. Mit seinen Mitarbeitern Jörg Voigt und Götz Krummheuer griff er dabei Lawler`s "Mikrowelt" auf und entwickelte sie weiter zum \emph{Subjektiven Erfahrungsbereich} (SEB) [Bauersfeld 1983].  

Während die Mikrowelt ein rein kognitives Modell im Sinne des "human information processing" ist, will das SEB -- Modell auch Erfahrungen erfassen und nicht nur Wissen. Damit werden auch prozedurales Wissen und nichtkognitive Dimensionen wie Motorik, Emotionen, Wertungen usw. mit eingeschlossen. \medskip  

Die Bereichspezifität der Konstruktion des mathematischen Wissens stützte sich nicht nur auf die Arbeit von Kognitionspsychologen sondern wurde u.a. auch durch die vielfältigen Bemühungen unterstrichen, bei denen der kulturelle Kontext besonders hervorgehoben wird
(vgl. [D`Ambrosio 1985, 1999], [Johnsen Hoines und Mellin--Olsen 1986]) --- zusammengefaßt unter dem Stichwort \emph{Ethnomathematik}.   \bigskip  

Eine für die Mathematikdidaktik zentrale Folgerung, die sich aus der Bereichspezifität oder Kontextgebundenheit der Wissenskonstruktionen ergibt, formuliert Bauersfeld wie folgt: 
\begin{quote}
\begin{small}
"Es gibt keine allgemeinen Begriffe. Zwar können wir in einer gegebenen Situation Begriffe mit jedem Allgemeinheitsgrad formulieren, aber sie werden damit nicht in gleicher Allgemeinheit im Gedächtnis abrufbar. Sie bleiben durch das Mitgelernte kontextgebunden, bereichsspezifisch" [Bauersfeld 1987] (vgl. auch [Bauersfeld 1985]).
\end{small}
\end{quote}  \medskip  

Als Bauersfeld seine Konzeption der SEBe veröffentlichte, reagierte die mathematikdidaktische Öffentlichkeit, soweit sie sich für Grundschulmathematik interessierte, überaus positiv. Miriam, der Name der Tochter Lawler`s, war damals nahezu jeder Grundschullehrerin bekannt. Auch heute noch liest sich der Artikel Bauersfelds faszinierend. Trotzdem ist seine Wirkung, rückwirkend betrachtet, eher begrenzt geblieben. Zwar werden von etlichen Autoren SEBe erwähnt, aber es wird stets auf die Lawlerschen Beispiele verwiesen. Wie aus dem oben Gesagten hervorgeht, hilft das, was einen SEB über eine Mikrowelt hinaus ausmacht, offenbar wenig, die SEB -- Konzeption zu einer präzise formulierbaren Theorie auszuarbeiten. Von daher stellt sich die Frage, ob es nicht erfolgversprechender ist, zunächst zu versuchen, Lawler`s Ausführungen über Mikrowelten (die "kognitiven Schatten" der SEBe [Bauersfeld 1983, S. 48]) zu präzisieren. Betrachtet  man seine charakterisierenden Beispiele (Miriams Geld -- Welt, Dekaden -- Welt, Papiersummen -- Welt und Zähl -- Welt), so genügen diese durchaus der oben angegebenen Charakterisierung der Zahlaspekte. \bigskip 

Eine Mikrowelt hat zwei Bestandteile, eine Perspektive und Funktionen. Die \emph{Perspektive} enthält Beschreibungen / Namen der Elemente, von denen die Mikrowelt handelt. Lawler trennt damit zwischen den Objekten der Mikrowelt und ihren Beschreibungen / Namen. Wir verzichten auf diese Trennung, da Miriam mit den Objekten der Mikrowelten agiert und nicht mit deren Beschreibungen / Namen. Wir sagen Beschreibungen statt Namen, da sich die Beschreibung durch das Arbeiten in der Mikrowelt zunehmend verfeinern kann. So unterschied Miriam in der Geld -- Welt zwei Sorten von Münzen zunächst nach ihrer Farbe und erst später nach ihrem Wert. Ein Name dagegen könnte sich nur ändern, kaum verfeinern. Vielleicht ist dies aber auch nur eine Frage des Sprachgebrauchs. Da eine verfeinerte Beschreibung die Objekte, von denen eine Mikrowelt handelt, nicht verändert --- eher etwas über die Verweildauer der Mikrowelt aussagt ---, verstehen wir die Perspektive einer Mikrowelt als die Menge der Objekte, über die sie Aussagen macht. \smallskip  
 
Die \emph{Funktionen} einer Mikrowelt sind (2-stellige) Relationen, definiert auf der Perspektive, oder Operationen, die sich auf die Elemente der Perspektive anwenden lassen. In Miriams Zählwelt sind dies z.B. die Relation = und die Operation +, angewandt auf die Zählzahlen, in der Geld -- Welt die Relation <. \medskip  

Eine \emph{Mikrowelt} T betrachten wir daher als ein (geordnetes) Paar (P,F) mit einer (nichtleeren) Menge P von Objekten und einer (nichtleeren) Menge F von (2-stelligen) Relationen oder Operationen, definiert auf P.  \medskip 

Unter diese Definition ordnen wir die (nach Bauersfeld) sog. \emph{Instrumentalen Welten} ein --- Lawler spricht von \emph{task -- rooted microworlds} ---, also die Zähl -- Welt, die Geld -- Welt, die Dekaden -- Welt und die Papiersummen -- Welt. Die Perspektiven dieser Welten "beschreiben (s.o.) Dinge unserer Alltagswelt" und "Vermutungen über deren Beziehungen" (die Funktionen, die Verf.) "können durch einfache Versuche im spezifischen Aufgabenfeld bestätigt oder verworfen werden" [ebd., S. 24].\vspace{0,5cm}

Bevor wir die weiteren Arten der Mikrowelten, Serien -- Welten und Konforme Welten, betrachten, greifen wir eine Einsicht auf, die Horst Struve in die Diskussion der Aspekte eingebracht hat. Einen Zahlaspekt zu beherrschen ist demnach damit gleichbedeutend, über eine empirische Theorie zu verfügen, deren Modelle die Klasse der Anwendungen ausmachen, bei denen die (natürlichen) Zahlen in gleicher Weise verwendet werden [Struve 1990, S. 202] (vgl. [Burscheid und Struve 1994]). Ersetzt man "Zahl" durch "Zahlzeichen", so läßt sich in gleicher Weise beschreiben, was es heißt, über eine Instrumentale Welt zu verfügen. Offenbar sind vergleichbare geistige Tätigkeiten erforderlich. In beiden Fällen erwirbt der Lernende eine empirische Theorie. Unter \emph{diesem Gesichtspunkt} macht es keinen wesentlichen Unterschied, ob man von Zahlaspekten oder von Instrumentalen Welten spricht. Mit dem Instrument "Empirische Theorie" lassen sich Zahlaspekte und Instrumentale Welten präzise beschreiben. \medskip

\noindent \emph{Bemerkung}: Wenn wir sagen, Kinder \emph{verfügen über} eine Theorie, so wollen wir damit ausdrücken, daß sie sich so verhalten, als hätten sie diese erworben, nicht, daß ihnen die Theorie in irgendeinem formalen Sinne bewußt wäre. Daß es angemessen ist, das Verhalten von Kindern so zu beschreiben, läßt sich schon bei Kleinkindern nachweisen (vgl. [Gopnik und Meltzoff 1997]). \medskip  

Das bereichsspezifische Wissen, das Kinder über Zahlaspekte oder Instrumentale Welten erwerben, ist in hohem Maße kontextgebunden und hat damit eine starke ontologische Bindung. Hält man den Wissenskorpus der Mathematik im Sinne Hilberts dagegen, so sind darin alle Begriffe von gleichem logischen Status, es sind Variable. Es ergibt sich zwangsläufig, daß bei einem bereichsspezifischen Wissenserwerb das Ergebnis der Wissenskonstruktion keine mathematischen Theorieteile sein können. Welches sind dann die Theorien, die Kinder erwerben?  \medskip

Eine kontextgebunde Auffassung der Wissenskonstruktion bei Kindern fragt zunächst nach diesem Kontext. Dieser ist in erster Linie gegeben durch die Umwelt, die die Kinder erfahren und bewältigen wollen und müssen. Die Wissenskonstruktion erstreckt sich daher zunächst auf den Entwurf von "Alltagstheorien", die sich im Umgang mit der WELT entwickeln und bewähren müssen. Wie gerade gesagt, verhalten sich schon Kleinkinder so, als verfügten sie über Theorien, die ihr Verhalten bestimmen. Alltagstheorien meint also solche Theorien (hier:) über mathematische Inhalte, die ihren Ausgangspunkt unmittelbar bei Problemen der Umwelt haben. Ihr Ziel ist es, ein Stück WELT zu beschreiben. Beispiele wären eine Theorie über das Verhalten eines Balles, eine elementare Meßtheorie oder eine Theorie über die gleichmäßige Zerlegung einer Schokoladentafel. In diesen Theorien haben die Begriffe nicht den Status von Variablen, sie haben vielmehr eine starke ontologische Bindung, und die Sätze sind keine Aussageformen, d.h. die Schülertheorien als mathematische Theorieteile darzustellen, wäre keine adäquate, auch den Lernprozeß in den Blick nehmende Beschreibung. Alltagstheorien entsprechen in ihrer Struktur und im Status ihrer Begriffe den Theorien von Naturwissenschaftlern. Wie einleitend schon gesagt, sprechen wir verallgemeinernd von \emph{empirischen Theorien}. Man kann sie als solche Theorien auffassen, die von realen Objekten (zumindest meßbaren Identitäten) handeln und deren Begriffe diese Objekte oder Bezüge zwischen diesen Objekten als Referenzen haben.  \medskip  

Es wurde schon darauf hingewiesen, daß das strukturalistische Theorienkonzept eine Möglichkeit bietet, solche Theorien formal zu beschreiben. In [Burscheid und Struve 2009] haben wir sowohl in allgemeiner Form die Beschreibung einer empirischen Theorie angegeben wie auch mehrere mathematikdidaktisch relevante Beispiele formal dargestellt, darunter diejenigen, die wir im folgenden heranziehen. Der wiederholt verwendete Terminus \emph{Theorie -- Element} weist darauf hin, daß empirische Theorien nur einzelne Anwendungen haben. Zusammenhängende Theorie -- Elemente werden zu \emph{Theoriennetzen} zusammengefaßt. Ein Beispiel bildet das Netz der Zahlaspekte, welches das Verfügen über den Begriff der natürlichen Zahl beschreibt.  \bigskip  
  
Den Erwerb der \textbf{Zählzahlen} kann man  als das Verfügen über ein Theorie -- Element (Z$;\preceq)$ beschreiben. Dabei bedeuten Z eine (nichtleere) Menge von Zahlwörtern (in konventioneller Reihenfolge) und $\preceq$ die (reflexive) lineare Ordnung der Zählzahlen, die durch Weiterzählen realisiert wird. \medskip  

Um Zählzahlen als \textbf{Maßzahlen} auffassen zu können wird auf Z$^3$ wie folgt eine Relation $\omega$ eingeführt:  

\begin{center}
$ \omega(n,m,l) \leftrightarrow $ Weiterzählen von n um m hat das Ergebnis l \medskip \\
\end{center}
Mit Hilfe von $ \omega $ läßt sich dann definieren:  
\begin{equation*}
n + m = l  \leftrightarrow  \omega(n,m,l) 
\end{equation*} 
   
\noindent (Z$; \preceq; +)$ bietet eine adäquate Beschreibung der \emph{Zähl -- Welt} als einer empirischen Theorie (vgl. [Burscheid und Struve 1994]). Ihre Grundbegriffe sind die Objekte und Funktionen der Lawlerschen Beschreibung der entsprechenden Mikrowelt.  \bigskip  

Da man in der Papiersummen -- Welt über die Zählzahlen hinaus die Zahl 0 ("Null") benötigt, erweist sich der Übergang von der Zähl -- Welt zur Papiersummen -- Welt als komplizierter als es auf den ersten Blick erscheinen mag. Um ihn korrekt darzustellen, müssen wir --- auch formal --- den Begriff der empirischen Theorie heranziehen (vgl. [Balzer 1982]). \bigskip \\
Wir definieren das Theorie -- Element \emph{Papiersummen} (PPS) wie folgt: \medskip 

\noindent M$_{pp}$(PPS)$\subseteq \langle \langle$Z;$ \preceq ;+\rangle \rangle $ , d.h. die \emph{partiellen Modelle} von PPS sind sämtlich Zähl  \medskip \\ -- Welten (genauer: Modelle des Theorie -- Elementes $\langle$ Z;$ \preceq ; + \rangle)$   \bigskip  \\
M$_p$(PPS) (die \emph{potentiellen Modelle} von PPS) $ \subseteq \langle\langle$ Z;$ \preceq ; + \rangle , 0 \rangle $ mit \bigskip  
\begin{enumerate}
 \item[(i)] $\langle$ Z;$ \preceq ;+ \rangle \in$ M$_{pp}$(PPS)  \smallskip  \\
\item [(ii)] Der linearen Ordnung (Z;$ \preceq $) wird ein kleinstes Element hinzugefügt, das 0 geschrieben wird und ein Zahlwort bezeichnen möge, das nicht unter den Zahlwörtern aus Z vorkommt. Statt Z$\, \cup \, \{0\}$ schreiben wir im folgenden $ \bar{\text{Z}} $. \smallskip 
\end{enumerate}
(0 wird somit als theoretisches Element in das Theorie -- Element PPS eingeführt.) \bigskip \\ 
M(PPS) (die \emph{Modelle} von PPS) = $ \langle \langle \langle$Z;$ \preceq ;+\rangle,0\rangle \rangle $ mit \medskip 
  \begin{enumerate}
\item [(i)] $ \langle \langle$Z;$ \preceq ;+ \rangle ,0 \rangle \in$ M$_p$(PPS)     \bigskip  
\item[(ii)] (Z;+) wird wie folgt ergänzt:  $ \underset{k \in \bar{Z}}{\bigwedge}(k + 0 = 0 + k = k) $  \bigskip  
 \end{enumerate}
 
Vergleicht man die Zähl -- Welt mit der Papiersummen -- Welt, so stellt man fest, daß die partiellen Modelle der letzteren Modelle der ersteren sind. Stehen zwei Theorie -- Elemente T und  T$^\prime$ in solch  einer Beziehung, so nennen Strukturalisten T$^\prime$ eine  \emph{Theoretisierung} von T, d.h. T$^\prime$ setzt T in dem Sinne voraus, daß jede intendierte Anwendung von T$^\prime$ unter Rückgriff auf die Begriffe und Gesetze von T beschrieben werden kann. Daß die Papiersummen -- Welt eine Theoretisierung der Zähl -- Welt dastellt,ist von ganz anderem Gewicht als die bei Lawler / Bauersfeld erwähnte Nachfolgerelation, zu deren Begründung lediglich darauf verwiesen wird, daß Miriam bei der noch zu behandelnden "Addition" in PPS die Spaltensummen durch Abzählen an ihren Fingern ermittelt [Bauersfeld 1983, S. 20]. Die genannte Theoretisierungsrelation, in der PPS und die Zähl -- Welt stehen, ist dagegen unabhängig von dem hier vorliegenden Beispiel und kann in jedem Theoriennetz zwischen einzelnen Theorie -- Elementen gegeben sein.  \bigskip  
  
Was die "Addition" in PPS betrifft, so läßt sich diese durch "paradigmatische Verhaltensbeschreibungen" [Rodenhausen 2010, S. 15] mittels geeigneter Tabellen einführen. Dabei wird spaltenweise addiert. Die Kenntnis der Operation + in $ \bar{Z} $ im Rahmen des "kleinen Eins + Eins" ist somit erforderlich (vgl.  M$_{pp}$(PPS)$\, \subseteq \, $M$ \big( \langle$ Z;$ \preceq ;+ \rangle \big) $. Die sog. Zehnerüberschreitung erfolgt rezeptartig, ohne Bezug zu einem Begriff wie Stellenwert o.ä.  \bigskip  
   
Als ein erstes Ergebnis können wir festhalten, daß auch zwischen Instrumentalen Welten, wenn man sie als empirische Theorie -- Elemente betrachtet, eine (nichttriviale) Nachfolgerelation bestehen kann. Dabei ist es durchaus möglich, daß verschiedene Instrumentale Welten einen gemeinsamen Vorgänger oder Nachfolger haben. \vspace{0,5cm}  

Als nächstes wenden wir uns dem Begriff der Serien -- Welt zu. Von Serien -- Welten wird gesagt, daß sie "wenigstens zwei Vorfahren" haben (In welchem Sinne?,  die Verf.) und --- was uns wichtiger erscheint ---, daß in ihnen "das unerwartete Zusammenwirken von Lösungswegen" deutlich wird [Bauersfeld 1983, S. 25]. Ohne daß es explizit gesagt wird, darf aus dem Kontext unterstellt werden, daß die erwähnten Vorfahren einer Serien -- Welt Instrumentale Welten sind, die durch sie zusammengefaßt werden. Da das erwähnte "Zusammenwirken von Lösungswegen" sich mathematisch im Begriff des Morphismus fassen läßt, definieren wir: \bigskip

T$_i$  = (P$_i$,F$_i$)  (i = 1,...,n) und T = (P,F) seien Instrumentale Welten  \medskip
 
T heißt \textbf{Serien -- Welt} für  T$_1$,...,T$_n \iff $ \bigskip 

$ \underset{i \in \mathbb{N}_n}{\bigwedge} \big[\underset{\alpha_i}{\bigvee} \underset{\beta_i}{\bigvee} \big((\alpha_i: P_i \rightarrow P,\: \text{injektiv,}\: \wedge \: \beta_i: F_i \rightarrow F,$ \text{injektiv und}  \bigskip 
 
\text{Stellenzahl neutral}$ \footnote{d.h. eine n -- stellige Relation wird auf eine n -- stellige Relation abgebildet}, \wedge \underset{a,b \in P}{\bigwedge} \underset{f\in F}{\bigwedge} \big((i) \vee (ii) \big) $ \medskip  
 
 (i) f Relation \bigskip 

 $ \big(a,b \in \alpha_i(P_i) \wedge f \in \beta_i(F_i) \wedge \big(\alpha_i^{-1} (a),\alpha_i^{-1}(b) \big ) \in \beta_i^{-1}(f) \rightarrow (a,b) \in f \big)\bigskip $ 
  
 (ii) f Operation \bigskip  
 
 $ \big(a,b \in \alpha_i(P_i) \wedge f \in \beta_i(F_i) \rightarrow \alpha_i\big(\beta_i^{-1}(f)(\alpha_i^{-1}(a),\alpha_i^{-1}(b)\big)\big) = f(a,b)\big) $ \bigskip
   
 Ist T eine Serien -- Welt für T$_1$,..., T$_n,$ so heißen T$_1$,..., T$_n$ \smallskip 
 
 \textbf{Spezialisierungen} von T
  \bigskip

\noindent 
 Gemäß dieser Definition ist die Zähl -- Welt eine Serien -- Welt für Geld -- Welt und Dekaden -- Welt. Geld -- Welt und Dekaden -- Welt sind Spezialisierungen der Zähl -- Welt. Als zweites Ergebnis erhalten wir somit, daß es keineswegs notwendig ist, der Serien -- Welt einen von der Instrumentalen Welt verschiedenen Status zuzuweisen, um die entscheidende Einsicht zu gewinnen, die eine Serien -- Welt vermitteln soll. \vspace{0,5 cm} 

Wir kommen zum Begriff "Konforme Welt". Auch die Konforme Welt verbindet Instrumentale Welten. Sie wird als  "relationale Welt" bezeichnet [ebd., S. 25]. Die aus unserer Sicht wesentliche Aussage, die über Konforme Welten gemacht wird, lautet: "Die Perspektiven der 'Konformen Welt' sind die Äquivalenzen zwischen den Perspektiven und Funktionen der beteiligten Mikrowelten" [ebd., S. 23]. Nimmt man das auf den Seiten 22 / 23 erläuterte Beispiel hinzu, so kann man die vorstehende Aussage auch wie folgt formulieren: Die Konforme Welt zeigt, daß Perspektiven und Funktionen der beteiligten Instrumentalen Welten nur unterschiedlich realisiert sind. Der Wissenszuwachs, den eine Konforme Welt vermittelt, ergibt sich aus der Zusammenschau der beteiligten Welten.

Dieses Problem läßt sich auch aus einem anderen Blickwinkel betrachten. Eine wesentliche Eigenschaft empirischer Theorien ist es, keine universellen sondern nur fest umrissene (Einzel--) Anwendungen zu haben. Natürlich hängen diese Anwendungen zusammen, denn sie gehorchen z.B. den gleichen Gesetzen. Dieser Zusammenhang wird im Begriffssystem empirischer Theorien in der Querverbindung ausgedrückt (vgl. Einleitung). Bezogen auf die hier aufgeworfene Frage nach einer darstellungsunabhängigen Problemlösung könnte sie beinhalten, daß es einer umfassenden Relation / Operation bedarf, die in jeder einzelnen Anwendung der dort vorkommenden Relation / Operation entspricht, d.h. es kommt nicht darauf an, wie die Relation / Operation in der einzelnen Anwendung modelliert wird.
Es scheint uns dies eine adäquate Sichtweise des Problems zu sein. Wir definieren. \bigskip  

T$_i$ =(P$_i$,F$_i)$ mit F$_i$ =$ \{f_{i\kappa} \rvert \kappa = 1,...,k\} $ (i = 1,...,n) seien Instrumentale Welten \medskip    

P = $\underset{i=1}{ \overset{n}{\bigcup}}$ P$_i$ \medskip
 
T$_1$,..., T$_n$ heißen  \textbf{konform modelliert} $ \iff $ \bigskip

$ \underset{\kappa \in \mathbb{N}_k}{\bigwedge} \underset {f_\kappa}{\bigvee} \big( f_\kappa \:  \text{Relation oder Opereration auf P} \wedge \underset{i \in \mathbb{N}_n}{\bigwedge} \big(f_{i\kappa} = f_\kappa/P_i \big) \big) $  \bigskip  

\noindent
 Wie es die Definition ausdrückt, betrachten wir die Konformität der Modellierung als eine Eigenschaft der Welten T$_1$,..., T$_n,$ nicht aber als eine eigene Mikrowelt. Das dritte Ergebnis lautet damit: Die Konforme Welt beschreibt ---  wie schon die Serien --  Welt ---  eine spezielle Eigenschaft der beteiligten Instrumentalen Welten. \medskip  

In Zähl -- Welt, Geld -- Welt. Dekaden -- Welt und Papiersummen -- Welt werden nur die Relationen $ \preceq $ und die Operation + betrachtet. Sie entsprechen in jeder dieser Welten der üblichen $ \leq $ -- Relation und der Addition natürlicher Zahlen. Damit sind die vier Welten konform modelliert.\bigskip  

Wie wir ausgeführt haben, sind die von Lawler angegebenen Instrumentalen Welten Modelle der Zählzahlen oder Spezialisierungen von ihnen bzw. --- im Falle der Papiersummen -- Welt --- eine Theoretisierung des zugehörigen Theorie -- Elementes, d.h. eigene Welten sind nur Instrumentale Welten. Um Lawler's Idee zu verifizieren, ist es nicht erforderlich, Serien -- und Konforme Welten als eigene Welten einzuführen. Sie lassen sich verstehen als spezielle Relationen zwischen Instrumentalen Welten.  \medskip 

Wir wollen noch kurz skizzieren, wie sich die Konzeption der Mikrowelten auf Brüche übertragen ließe. Bekannte Bruchzahlaspekte sind Maßzahlen, Operatoren und Verhältniszahlen. Zugehörige Instrumentale Welten wären z.B. Verteilungsprobleme (von Pizzen oder Schokoladentafeln), Maßstäbe von Plänen oder Landkarten und Einkaufssituationen oder auch spezielle Verteilungsprobleme (Verteilen von Anteilen).   

In [Burscheid und Struve 2009] haben wir das Theorie -- Element einer \emph{Empirischen Meßtheorie mit rationaler Skala} angegeben. In ihm wurde eine Multiplikation eingeführt, die wir als \emph{verkettete Teil -- Ganzes -- Beziehung} bezeichnet haben. (Auf eine detaillierte Angabe verzichten wir hier wegen des erforderlichen Aufwandes.) Betrachtet man ein sog. \emph{inneres Verhältnis} zweier Objekte, das Verhältnis einer übereinstimmenden Eigenschaft dieser Objekte, und erfaßt man die Bestimmung dieses Verhältnisses durch eine empirische Theorie, so läßt sich zeigen, daß die Modelle des zugehörigen Theorie -- Elementes auch Modelle des Theorie -- Elementes der \emph{Empirischen Meßtheorie mit rationaler Skala und verketteter Teil -- Ganzes -- Beziehung} sind, d.h. das Verständnis von und der Umgang mit inneren Verhältnissen läßt sich innerhalb der \emph{Empirischen Meßtheorie} erwerben. Wie die \emph{Lernsequenz über Bruchzahlen als Maßzahlen} zeigt, in die die \emph{Empirische Meßtheorie} umgesetzt wurde, läßt sich auch der Umgang mit Operatoren (über das Verteilen von Anteilen) innerhalb der \emph{Empirischen Meßtheorie} entwickeln. Beide Theorie -- Elemente, dasjenige über den Umgang mit inneren Verhältnissen und dasjenige über den Umgang mit Operatoren sind --- was sich formal nachweisen läßt --- Spezialisierungen der Instrumentalen Welt der Verteilungsprobleme oder --- anders gewendet --- die Instrumentale Welt der Verteilungsprobleme ist eine Serien -- Welt der Instrumentalen Welt der Maßstabsbestimmungen und derjenigen der Verteilung von Anteilen. \medskip

Die angegebene Präzisierung des Lawlerschen Konzeptes ist also nicht auf auf die Analyse eines ersten Zahlverständnisses beschränkt, nicht einmal auf die natürlichen Zahlen, sondern kann in systematischer Weise genutzt werden, bei unterschiedlichen Zahlbegriffen Zusammenhänge zu untersuchen. Wie wurde dies erreicht? Einmal dadurch, den Erwerb einer Instrumentalen Welt als den Erwerb einer empirischen Theorie aufzufassen, zum andern dadurch, die Konzepte Serien -- Welt und Konforme Welt nicht von ihrer Definition sondern \emph{von dem Zweck her} zu bestimmen, dem sie dienen sollen. Daraus ergibt sich die Möglichkeit, sie als intertheoretische Relationen aufzufassen, eine deutliche Abweichung vom Lawlerschen Original. Im nachhinein sieht man jetzt auch den Bruch im Lawlerschen Original. Die Instrumentalen Welten haben einen empirischen Anteil, Serien -- Welten und Konforme Welten aber sind mentale Modelle. Wie man sieht, kann man auf die letzteren verzichten, will man die Entwicklung der Instrumentalen Welten beschreiben.

Die von Poser genannte Gefahr, das Selbstverständnis des Konzeptes zu rekonstruieren, besteht deshalb nicht, weil kein alternatives Vorverständnis des Lawlerschen Konzeptes vorliegt. Zumindest ist uns keines bekannt. Unser vornehmliches Interesse war es, das Konzept --- mit den Worten von Poser --- für die Behandlung \emph{systematischer} Probleme dienstbar zu machen.

Zur Frage der Adäquatheit der Rekonstruktion ist zu sagen, daß die Antwort natürlich davon abhängt, was man als Rekonstruktion gelten läßt. Für jemanden, der nur Direktinterpretationen akzeptiert, ist die von uns vorgenommene Rekonstruktion natürlich nicht adäquat. Wer der Verwendbarkeit des Konzeptes der Vorzug gibt, dürfte die Adäquatheit der Rekonstruktion als gegeben ansehen, da keine Argumente eingeflossen sind, die dem Lawlerschen Konzept nicht angehören. \medskip  

Daß unseres Wissens nach dem Erscheinen der Lawlerschen Arbeit niemand neue Mikrowelten in systematischer Weise benannt hat, dürfte damit zusammenhängen, daß sie --- da äquivalent zu Zahlaspekten --- nur in diesem Kontext auftreten. Die Diskussion der verschiedenen Aspekte natürlicher Zahlen war aber schon weitgehend abgeschlossen als die Mikrowelten durch Bauersfeld bekannt gemacht wurden. Ein Versuch, Lawler`s Beobachtung bei anderen Zahlbegriffen als dem der natürlichen Zahlen zu untersuchen, ist uns nicht bekannt.
\vspace{0.5cm}  

Die von Bauersfeld verfolgte Intention, eine "Interaktionstheorie mathematischen Lernens und Lehrens" zu skizzieren [ebd., S. 27], wird durch diesen Beitrag nicht berührt, da eine freizügige Verwendung des Begriffs SEB dazu nicht erforderlich ist. Gegen diese spricht dann nichts, wenn man SEBe nur als mentale Modelle betrachtet. Dies ist aber nicht unproblematisch, da man vorgibt, von dem empirischen Gehalt der Instrumentalen Welten absehen zu können.

\section{\textbf{Die Einarbeitung empirischer Befunde}}
Wie wir schon in der Einleitung angedeutet haben, ist es sehr unbefriedigend, daß empirische Befunde nur sehr indirekt Einfluß auf die Gestaltung der Curricula nehmen. Und wie ebenfalls schon gesagt wurde, sehen wir eine Möglichkeit, diesem Umstand abzuhelfen, darin, den Wissenserwerb des Schülers --- im systematischen Sinne --- durch eine empirische Theorie zu beschreiben. Dazu genügt es allerdings nicht, nur deren begriffliche Struktur zu nutzen wie bei dem vorstehenden Beispiel, sondern es muß auch die formale Darstellung empirischer Theorien herangezogen werden, um jeden erforderlichen Schritt detailliert angeben zu können. Da der Erwerb des Wissens nachgezeichnet und nicht nur das Ergebnis in mathematischer Form festgehalten wird, können empirische Befunde berücksichtigt werden. Der Präzisionsgrad einer solchen Darstellung erlaubt dann eine rational geführte Diskussion über unterrichtliche Konsequenzen.

Das Beispiel, an dem wir das Gesagte exemplifizieren wollen, ist der Maßzahlbegriff für ebene Figuren. Liest man empirische Untersuchungen zum Verständnis des Begriffs Flächeninhalt oder zur Bestimmung des Flächeninhaltes gegebener Figuren, so weisen die Autoren in der Regel bei den Probanden beachtliche Defizite nach. Dies vom Grundschüler [Hirstein u. a. 1978] bis zum angehenden Primarstufenlehrer [Tierney u.a. 1990]. Obwohl diese Defizite empirisch erhärtet sind und immer wieder beklagt werden, muß man fragen, inwieweit Unterrichtskonzeptionen ihnen Rechnung tragen. Im folgenden werden wir ein Theorie -- Element angeben, das die Bestimmung von Flächeninhalten zum Gegenstand hat und etliche der empirisch nachgewiesenen Defizite aufgreift.  \medskip

Die Schwierigkeiten des Themas Flächeninhalt beginnen schon im Sprachlichen. Während man von der Fläche eines Gartens spricht und diese durch den den Garten umgebenden Zaun begrenzt wird, der Zaun also den Rand der Gartenfläche ausmacht --- Ähnliches gilt für eine bebaute Fläche, die Wasserfläche eines Sees oder die Wohnfläche eines Hauses ---, ist es bei geometrischen Figuren leider so, daß umgangssprachlich nicht konsequent zwischen Fläche und Rand einer Fläche unterschieden wird. "Dreieck" kann sowohl drei nichtkollineare Punkte nebst den sie verbindenden Strecken meinen wie auch Dreiecksfläche bedeuten. Wir sprechen vom "Inhalt eines Dreiecks", seltener von der "Größe der Dreiecksfläche". Um solche Mißverständnisse auszuschließen, legen wir fest: Wenn wir die Fläche des Dreiecks meinen, sprechen wir von der \emph{Dreiecksfläche}. Die Strecken, welche die drei nichtkollinearen Punkte verbinden, bezeichnen wir als \emph{Dreieck} oder \emph{Rand der Dreiecksfläche}. Entsprechernd verfahren wir bei anderen geometrischen Figuren. Mit \emph{Figuren} sind einfach geschlossene ebene Polygone gemeint wie Rechtecke, Dreiecke, Parallelogramme o.ä. 

Um eine Theorie anzugeben, die die Bestimmung von Flächeninhalten zum Gegenstand hat, formulieren wir in einem ersten Schritt eine (Vor --) Theorie, deren Beherrschung das Wissen umfaßt, über das Kinder verfügen müssen, wollen sie sich erfolgreich mit der Bestimmung von Flächeninhalten beschäftigen. Sie soll vor allem das begriffliche Verständnis der später zu entwickelnden Berechnungsformeln (für Rechtecke, Dreiecke, Parallelogramme, ...) fördern [Huang und Witz 2011].

Bevor wir uns der Darstellung der Theorie zuwenden, weisen wir darauf hin, daß die Begriffe \emph{Fläche} und \emph{Rand einer Fläche} nicht im Rahmen der Theorie erworben werden, sondern es wird vorausgesetzt, daß die Kinder sich diese schon bei der Auseinandersetzung mit ihrer Umwelt --- an Beispielen wie den oben genannten --- angeeignet haben und zur Beantwortung von Fragen benutzen, die ihr Umfeld ihnen stellt. So bekommen sie Verständnis für die Zweidimensionalität von Flächen durch das Bestimmen geometrischer Örter. Um einen Punkt auf einer Fläche festzulegen, benötigen sie ein Maß auf zwei (vorgegebenen) orthogonalen Geraden. Die beiden Maße sicher zu koordinieren gelingt ihnen in der Regel erst ab dem 9. Lebensjahr [Piaget u. a. 1974]. So entwickelt sich die Vorstellung eines Netzes sich rechtwinklig schneidender Geraden, mit der man eine Fläche überziehen kann. Darin drückt sich für Kinder die Homogenität der Fläche aus. Flächen werden so beliebig teilbar, sind nicht unzerstörbare diskrete Objekte [Kamii und Kysh 2006].

Die Notwendigkeit, das Bestimmen von Flächeninhalten auf enaktiver Stufe vorzubereiten, betonen Thomas P. Carpenter --- er spricht von \emph{"premeasurement, establishing empirical procedures for directly composing, ordering, and combining elements of some domain of elements that possess a given attribute; ... "} (Carpenter o.J. [1976]) --- und Alan R. Osborne. Bei ihm heißt es: \emph{"Several different attributes of the domain set"} (der Flächeninhaltsfunktion, die Verf.) \emph{"may be operationally defined. Among these are:} 
\begin{quote}
\begin{small}
1. \emph{Transitivity property. If region A has the same area as region B and region B has the same area as region C, then region A has the same area as region C} (\emph{similarly, for less area than and greater area than}).\\ 
2.\emph{Substitutive property. If region A has the same area as region B and region B has more} (\emph{less}) \emph{area than region C, than region A has more} (\emph{less}) \emph{area than region C.} \\
3. \emph{Symmetric property. If region A has the same area as region B, then region B has the same area as region A.} \\
4. \emph{Asymmetric property. If region A has more} (\emph{less}) \emph{area than region B, then region B does have not more} (\emph{less}) \emph{area than region A}. \\ 
5. \emph{Reciprocal property. If region A has less area than region B, then region B has more area than region A."} (Osborne o.J. [1976]) 
\end{small}
\end{quote} 
Da die intendierten Anwendungen der Theorie --- die Objekte, über die die Theorie Aussagen macht --- empirische Objekte sind, Realisanten derjenigen Figuren, deren Flächeninhalte von Kindern als erste bestimmt werden --- Rechtecke, Quadrate, Dreiecke o.ä. --- und ihre Aussagen Bezüge zwischen diesen Objekten thematisieren, ist die Theorie, die wir  angeben werden, eine \emph{empirische} im oben genannten Sinne und \emph{nicht} eine mathematische im heutigen Verständnis. Wir sprechen von der \emph{Empirischen Theorie des Flächenvergleichs}. \\  

Wie wir in der wiederholt genannten Monographie näher ausgeführt haben, sollte jede didaktische Bearbeitung mathematischer Inhalte bestimmte \emph{methodologische Forderungen} erfüllen. Sie lauten: \smallskip

--- Es ist der Zweck anzugeben, zu dem der Schüler den neu einzuführenden Inhalt erlernen soll;

--- es ist anzugeben, wie der neu einzuführende Inhalt dem genannten Zweck dienen kann; 

--- es sind die systematischen Voraussetzungen detailliert anzugeben, auf die zurückgegriffen wird. \\

Die Antworten auf die methodologischen Forderungen nennen wir \emph{didaktische Postulate}. Ihnen das zu formulierende Theorie -- Element --- die \emph{Empirische Theorie des Flächenvergleichs} --- zu genügen. \emph{Didaktische Postulate}, da sie die Grundlage einer didaktischen Konzeption bilden, die präskriptive Aussagen beinhaltet.  \smallskip

(P1)  Der Zweck, zu dem ein Schüler den Vergleich von Flächengrößen erlernt, ist die Bewältigung von Alltagsproblemen.

(P2)  Die Handhabung qualitativer Vergleiche und Verfahren, Flächen auszulegen --- mit greifbaren Objekten und vor allem durch Anfertigung von Zeichnungen --- und über ihre Gesetzmäßigkeiten zu verfügen, ist eine notwendige Voraussetzung für eine später zu erwerbende  \emph{Empirische Theorie der Flächeninhaltsbestimmung}.
 
(P3)  Systematische Voraussetzungen für eine \emph{Empirische Theorie des Flächenvergleichs} sind 
 
--- die Fähigkeit, Gebilde (s.u.) auch über die Grenzen verschiedener Anwendungen hinweg als dieselben zu erkennen,
 
--- die Begriffe Fläche und Rand einer Fläche --- wie oben gesagt --- auf empirischem Wege erworben zu haben. \medskip
 
Die unterrichtliche Behandlung verschiedener wichtiger Inhalte stützt sich auf standardisierte Hilfsmittel, so die Behandlung des Zählens auf Perlen und Plättchen, die Einführung der Bruchrechnung auf Schokoladentafeln oder Pizzen usw. Die \emph{Empirische Theorie des Flächenvergleichs} beschreibt für die hier anstehende Thematik ein uns geeignet erscheinendes Hilfs -- oder Arbeitsmittel sowie die Einsichten, die es vermittelt. Diese Einsichten bilden auf enaktiver Stufe eine Basis für die anschließend beschriebene \emph{Empirische Theorie der Flächeninhaltsbestimmung}. \smallskip
 
Wir beschreiben die \emph{Empirische Theorie des Flächenvergleichs} durch ein vollständig definiertes Theorie -- Element. Wir bezeichnen es mit $T_{FV}$ (FV für "Flächenvergleich"). 
 
Wir beginnen die Definition des Theorie -- Elementes $T_{FV}$ mit der Angabe der partiellen Modelle. Dazu treffen wir folgende Vereinbarung:  

\noindent u$_i$ $(i \in \mathbb{N})$ seien eine hinreichend große Zahl  dreieckiger, rechteckiger, quadratischer, parallelogrammförmiger, ... \emph{Bausteine}: Plättchen, Kartonkärtchen, Bauklötze (geringer Dicke) o.ä. (vgl. [Doig et al. 1995]).

\noindent Für k $\in \mathbb{N}$ bezeichne [u$_1$,...,u$_k]$ ein Muster, das man erhält, indem man die u$_\kappa$ auf einer ebenen Fläche so ausbreitet, daß sie sich weder ganz noch teilweise überdecken. Die Muster [u$_1$,...,u$_k]$ nennen wir \emph{Gebilde}, um die Bezeichnung "Figur" den geometrischen Figuren vorzubehalten. Gebilde bezeichnen wir mit U, V o.ä. \smallskip

\noindent Die \emph{partiellen Modelle} des Theorie -- Elementes: \\

M$_{pp}$(T$_{FV}) = \big \langle \langle \Delta; \; \cong, \; \subsetneqq; \; \circ \rangle \big \rangle $ mit \smallskip 
\begin{enumerate}
\item [(i)]$ \Delta $: eine Menge von Gebilden, hergestellt aus Bausteinen u$_1$,...,u$_k$ gemäß vorstehender Vereinbarung \\ 
\item[(ii)]$\cong \; \subseteq \Delta \times \Delta$, aufgefaßt als "\emph{deckungsgleich}" im Sinne des Aufeinanderlegens \\
\item[(iii)]$\circ \subseteq \Delta^3$, aufgefaßt als "\emph{ergänzen eines Gebildes um ein zweites}"\footnote{"\emph{ergänzen}" heißt "hinzufügen des zweiten Gebildes zu dem ersten, ohne daß sich Bausteine ganz oder teilweise überdecken".} \\
\item [(iv)] $ \subsetneqq $, aufgefaßt als "\emph{ist echtes Teilgebilde von}", d.h. "$ U \subsetneqq V \leftrightarrow $ U läßt sich durch Hinzufügen von Bausteinen zu V ergänzen", formal: durch wiederholtes Anwenden von $\circ$ \\
\end{enumerate}

\noindent \emph{Bemerkung}: \\ 
\noindent 1. Die Bausteine werden wie abgeschlossene geometrische Figuren behandelt. \\
2. Um die Schreibweise nicht unnötig zu komplizieren, benutzen wir die üblichen mathematischen Symbole, auch wenn sie sich auf Gebilde oder deren ikonische Realisationen (durch Zeichnungen) beziehen. \\
3. Unter Rückgriff auf $\circ$ läßt sich die Konstruktion eines Gebildes [u$_1$,...,u$_k]$ wie folgt fassen: 
\begin{center}
[u$_1$, ...,u$_k]$ =$\big[$...[[u$_1$] $\circ$ u$_2$]...$\big]\, \circ$ u$_k$ $ \big]$
\end{center}
\noindent Ist U ein Gebilde, sind v$_1$,...,v$_l$ Bausteine, so ist folglich [U,v$_1$,...,v$_l$] ebenfalls ein Gebilde. \\ \bigskip

In den partiellen Modellen läßt sich wie folgt eine Relation zgl $\subseteq \Delta \times \Delta$ ("\emph{ist zerlegungsgleich zu}") definieren: \medskip \\
U zgl V $\leftrightarrow \underset{k \in \mathbb{N}}{\bigvee} \; \underset{u_1, ...   u_k \in \Delta}{\bigvee} \; \underset{\alpha}{\bigvee}(\alpha: \mathbb{N}_k \rightarrow \mathbb{N}_k, \text{bijektiv},\smallskip \\  \noindent \hspace*{11mm} \wedge \; U \cong [u_1, ... ,u_k] \wedge V \cong [u_{\alpha(1)}, ... ,u_{\alpha(k)}])$  \smallskip \\ 
\emph{Bemerkung}: Für k $\in \mathbb{N}$ bezeichne $ \mathbb{N}_k$ die Menge der Zahlen 1,...,k. \bigskip \\ 
\noindent Die \emph{potentiellen Modelle}: \\ 

M$_p$(T$_{FV})$ =$ \big \langle \langle \Delta; \; \cong, \; \subsetneqq, \; \text{zgl}, \prec; \; \circ \rangle \big \rangle$ mit \\
\begin{enumerate}
\item [(i)] $\langle \Delta; \; \cong, \; \subsetneqq, \; \text{zgl}; \; \circ \rangle \in$ M$_{pp}$(T$_{FV})$  \\ 
\item [(ii)] Seien U, V $\in \Delta$ \medskip \\
 Def.: U \emph{hat einen kleineren Flächeninhalt als} V $ \leftrightarrow \underset{V'\in \Delta}{\bigvee}(U \, \text{zgl} \, V' \wedge V' \subsetneqq V)$ \\
\noindent \hspace*{8mm} bez.: U $\prec$ V \medskip \\
\noindent \hspace*{8mm} U \emph{hat den gleichen Flächeninhalt wie} V $\leftrightarrow$ U zgl V \bigskip
\end{enumerate}
\noindent Und die \emph{Modelle}: \\

M(T$_{FV}) = \big \langle \langle \Delta, \; \cong, \; \subsetneqq, \; \text{zgl}, \; \prec; \; \circ \rangle \big \rangle$ mit \\
\begin{enumerate}
\item [(i)] $\langle \Delta; \; \cong, \; \subsetneqq, \; \text{zgl}, \; \prec; \; \circ \rangle \in$ M$_p$(T$_{FV})$ \smallskip \\
\item[(ii)] $\cong $ ist reflexiv, symmetrisch und transitiv auf $\Delta$ \\
\item[(iii)] $\subsetneqq$ ist asymmetrisch und transitiv (damit irreflexiv) auf $\Delta$ \\
\item [(iv)] zgl ist reflexiv, symmetrisch und transitiv auf $\Delta$ \\
\item[(v)] $\underset{U,V,W \in \Delta}{\bigwedge} \big (U \circ(V \circ W)\; \text{zgl} \; (U \circ V) \circ W \big )$ \\ 
\item[(vi)] $ \underset{U,V \in \Delta}{\bigwedge}(U \circ V \; \text{zgl} \; V \circ U)$ \\
\item[(vii)] $ \underset{U,V \in \Delta}{\bigwedge} (U \prec V \vee U \, \text{zgl} \, V \vee V \prec U))$ \smallskip \\
\emph{Bemerkung}: Wegen zgl $\, \cap \, \subsetneqq \: = \: \emptyset$ schließen die Bedingungen einander aus.  \\
\item[(viii)] $ \Delta_\square $ bezeichne die Menge aller rechteckigen Gebilde aus $\Delta$ \smallskip \\
$ \underset{U \in \Delta}{\bigwedge} \; \underset{V \in \Delta}{\bigvee} \;$ (U zgl V $\wedge$ V $\in \Delta_\square) $ \\
\end{enumerate} \medskip

 Als unmittelbare Folgerung erhält man \medskip \\
(T$_{FV}:1)$ Beh.: $\,\prec \,$ ist asymmetrisch (damit irreflexiv) und transitiv auf $\Delta$ \bigskip

Wir fragen nach der Adäquatheit der Modellaxiome. Die Axiome (ii) - (vi) sind so zu verstehen, daß die Schüler die jeweilige Relation so handhaben, als hätte sie die notierten Eigenschaften, nicht aber, daß  diese für sie in einem formalen Sinne überprüfbar wären. (vii) und (viii) sind die entscheidenden Axiome. Flächeninhalt ist ein bezüglich des Theorie -- Elementes T$_{FV}$ theoretischer Begriff. Er ist ein komparativer Begriff und (vii) verlangt die durchgängige Vergleichbarkeit der Gebilde hinsichtlich ihres Flächeninhaltes. Nur empirische Untersuchungen können zeigen, ob die Forderung zum Theorie -- Element T$_{FV}$ "paßt". (viii) dient zur Vorbereitung der \emph{Flächenformel} (Inhaltsformel des Rechtecks), deren breite Anwendbarkeit die Verwandlung einer geometrischen Figur in ein Rechteck voraussetzt. Die Wichtigkeit von (vii) betonen Osborne (Osborne o.J. [1976]), Constance Kamii und Judith Kysh  [Kamii und Kysh 2006] sowie James J. Hirstein, Charles E. Lamb und Osborne, die die aus der Zerlegungsgleichheit folgende Inhaltsgleichheit als den entscheidenden Schritt zum Verständnis des Begriffs Flächeninhalt ansehen. Sie sind der Auffassung, daß das Zerlegen von Flächen und das Zusammensetzen der entstehenden Teile zu anderen Formen beherrscht werden müsse wie die elementaren Rechenoperationen [Hirstein u. a. 1978]. Darüberhinaus ist (vii) eine Reichhaltigkeitsforderung an das gewählte Arbeitsmittel.  \medskip

Um die formale Darstellung des Theorie -- Elementes T$_{FV}$ abzuschließen bedarf es noch der \emph{Querverbindung}. Sie wäre entsprechen dem ersten Beispiel --- bezogen auf die Relationen $\cong, \, \subsetneqq, \,$ zgl und $\circ$ --- zu formulieren. Auf ihre Angabe werden wir auch hier verzichten. \bigskip

Man kann der Frage nachgehen, wie sich das formulierte Theorie -- Element von den (mathematisch fundierten, wenngleich genuin) physikalischen Theorien des Messens (von Längen und Flächen) unterscheidet, die sich ebenfalls als empirisch verstehen.  Als Beispiel einer solchen Theorie betrachten wir diejenige, die in [Krantz u. a. 1971] als "closed extensive structure" bezeichnet wird. Sie wird dort wie folgt definiert: \\

\noindent $\langle$A;\,$\succeq\, ; \circ \rangle$ heißt \emph{closed extensive structure}, wenn folgende Axiome für alle a, b, c, d $\in$ A gelten
\begin{enumerate}
\item [1.] $\succeq$ ist reflexiv, transitiv und konnex \smallskip \\
\emph{Bemerkung}: Definiert man $a \sim b \leftrightarrow a \succeq b \; \wedge \; b \succeq a \; \text{und} \; a \succ b \leftrightarrow a \succeq b \; \wedge \; b \nsucceq a$, so ist $\sim$ eine Äquivalenzrelation und $\succ$ eine irreflexive lineare Ordnung auf A \\
\item[2.] $\circ$ ist eine Operation auf A mit $a \circ (b \circ c) \sim (a \circ b) \circ c$ \\
\item[3.] $a \succeq b \leftrightarrow a \circ c \succeq b \circ c \leftrightarrow c \circ a \succeq c \circ b$ \\
\item [4.] $ a \succ b \rightarrow \underset{c,d \in A}{\bigwedge} \; \underset{n \in \mathbb{N}}{\bigvee} (na \circ c \succ nb \succ d)$ \\
(mit 1a = a und (n + 1)a = na $\circ$ a) \\
\item[5.] $(a \circ b \succ a)$ \\ 
\end{enumerate} 
 \noindent Folgende Interpretation wird dazu angegeben: 
 
 A sei eine Menge von Objekten, die eine bestimmte Eigenschaft besitzen. a $\succeq$ b bedeute, daß die Eigenschaft bei a mindestens so stark ausgeprägt ist wie bei b. a $\circ$ b ist ein Objekt aus A, das aus a und b in vorgeschrieben geordneter Weise zusammengefügt wird.
 
 Als Beispiel wird das Wiegen mit einer Balkenwaage genannt: A sei eine Menge von Objekten, deren jedes die Eigenschaft hat, schwer zu sein, und a $\circ$ b bedeutet, daß a und b in derselben Schale --- a unter b --- liegen.\medskip
 
 Was den empirischen Charakter dieser Struktur betrifft, so sehen wir keinen entscheidenden Unterschied zu unserem Ansatz. Hier hat man Objekte mit der Eigenschaft, schwer zu sein, dort hat man Gebilde mit der Eigenschaft, eine bestimmte Fläche auszufüllen. Hier hat man das Zusammenlegen der Objekte auf eine Schale, dort hat man das Ergänzen der Gebilde. Ein wesentlicher Unterschied auf der empirischen Ebene ist für uns nicht erkennbar. Die Abweichungen der beiden Axiomengruppen, die die empirischen Vorgaben präzisieren, sind durch die verschiedenen Inhalte bedingt bzw. Konzessionen an das Alter der Schüler (z.B. $\subsetneqq$ statt $\subseteqq$). 
 
 Ein deutlicher Unterschied liegt allerdings in der Zielsetzung der beiden Theorien. Das Ziel von Krantz u. a. ist es, eine mathematische Struktur zu beschreiben, die physikalische Systeme erfüllen müssen, damit man Messungen mit reellen Zahlen vornehmen kann, genauer: damit ein Homomorphismus von A in $\mathbb{R}$ existiert. Dazu definiert er eine solche Struktur (im modernen Sinne).
 
 Unser Ziel ist es, einen Lernprozeß (im systematischen Sinne) zu beschreiben. Dies wird deutlich \\
 - im Stufenaufbau des Theorie -- Elementes, insbesondere der Unterscheidung von theoretischen und nichttheoretischen Begriffen und der Identifizierung von Lernschritten, \\
 - im Interesse an den Referenzbeziehungen der verwandten Begriffe und der Überprüfbarkeit der Axiome (Setzungen durch Axiome vs. empirisch überprüfbare Aussagen), \\
 - in der Rechtfertigung der Axiome durch empirische Untersuchungen (und nicht durch mathematische Überlegungen; die Rechtfertigung der \emph{closed extensive structure} besteht im Nachweis ihrer Darstellbarkeit mit Hilfe reeller Zahlen), \\
 - in der Spezifizierung der Voraussetzungen für den Erwerb der Theorie.\medskip
 
 Die unterschiedlichen Ziele werden deutlich in der Formulierung der Axiome. Währen sie für die \emph{closed extensive structure} möglichst elegant, sparsam und in einem mathematischen Sinne vollständig sein sollen, dienen die Modellaxiome des Theorie -- Elementes T$_{FV}$ nur dazu, den empirischen Vorgaben einen präzisen begrifflichen Rahmen zu geben. Die Aufgabe der Modellaxiome ist lediglich, \emph{im Rahmen des Kontextes} geeignete Modelle einzugrenzen, um die Zuweisung von Größenwerten in bekannter Weise vornehmen zu können. Deshalb wird nicht nach der Unabhängigkeit der Axiome oder der Minimalität des Systems gefragt, wie es in einer mathematischen Theorie üblich ist. \bigskip
 
 Was leistet ein Arbeitsmittel, welches das Theorie -- Element T$_{FV}$ modelliert? Genauer: Was soll es leisten? Zunächst sei auf die oben von Osborne formulierten Eigenschaften verwiesen, die die Definitionsmenge der Flächeninhaltsfunktion --- im vorliegenden Falle beschränkt auf Polygonflächen --- besitzt. Diese seien operational zu definieren, was heißt, daß sie sich dem Schüler --- sich wechselseitig bedingend --- in enaktivem Vollzug erschließen sollen. Wie man an Hand der Modellaxiome verifiziert, beinhaltet, über das Theorie -- Element T$_{FV}$ zu verfügen, die Einsichten gewonnen zu haben, die Piaget der premeasurement -- Stufe zuweist.

Hervorgehoben werden sollte noch, daß die Transitivität (im Sinne Piagets "qualitativer Transitivität" [Piaget u. a. 1974]) und die "conservation ability" von Flächengrößen --- letztere in der Form, daß bestimmte Änderungen, denen ein Gebilde unterworfen wird, sodaß es eine andere Gestalt erhält, z.B. indem man es zerlegt und in neuer Weise zusammensetzt, seinen Flächeninhalt invariant läßt --- sehr nahe beieinander liegen (Carpenter o.J. [1976]) und das Verständnis der einen Eigenschaft das der anderen fördert [Taloumis 1975]. Darauf hinzuweisen ist deshalb wichtig, weil die conservation ability wesentlich ist, um zu verstehen, daß nichtkongruente Gebilde inhaltsgleich sein können. Eine Einsicht, über die nicht alle Schüler verfügen [Hirstein u. a. 1978], [Kamii und Kysh 2006]. \medskip

\noindent \emph{Bemerkung}: Daß die conservation ability für das Verständnis der Flächeninhaltsbestimmung äußerst wichtig ist, betonen auch Harriett G. Wagman [Wagman 1975] und Jacques Lautrey, Etienne Mullet, Patricia Paques [Lautrey u. a. 1989] sowie in einer neueren Untersuchung George Kospantaris, Panagotios Spyrou, Dionyssios Lappas [Kospantaris u. a. 2011]. \medskip

Da es ein Ziel der Theorie ist, die Zerlegungsgleichheit von Gebilden als Gleichheit von Flächengrößen zu verstehen, impliziert die Transitivität von zgl das "Substitutionsprinzip" in folgender Form: Man hat zwei Gebilde gleichen Flächeninhaltes. Vergleicht man ein drittes mit diesen hinsichtlich der Größe des Flächeninhaltes, so erhält man jeweils das gleiche Ergebnis.\medskip

Weitere Einsichten, die die Theorie vermittelt, sind

--- die Teil -- Ganzes -- Beziehung: Der Flächeninhalt eines Teilgebildes ist kleiner als der Flächeninhalt des Gesamtgebildes (vgl. M(T$_{FV}$) (vii));

--- das Verständnis der sich bedingenden Beziehungen "Gebilde U hat einen kleineren Flächeninhalt als Gebilde V" g.d.w. "Gebilde V hat einen größeren Flächeninhalt als Gebilde U" (vgl. M(T$_{FV}$) (iii) und (iv)). \medskip

Empirische Untersuchungen zum Thema Flächeninhalt heben hervor, daß auch das Auslegen z.B. einer Rechteckfläche mit Einheitsquadraten beim Schüler nicht notwendig die Vorstellung hervorruft, daß dabei eine Fläche ausgefüllt werde. Die Einheitsquadrate werden als diskrete Ganzheiten betrachtet, die man zählen aber z.B. nicht teilen kann. (Harry Beilin (1964) nach Osborne o.J. [1976]), [Kamii und Kysh 2006], [Hirstein u. a. 1978]. Ein T$_{FV}$ modellierendes Arbeitsmittel versucht dem entgegenzuwirken, indem die Grundmenge $\Delta$  unterschiedlich gestaltete Bausteine (Dreiecke, Rechtecke, ...) enthält.

Dieses Problem dürfte auch stark von der verwendeten Begrifflichkeit abhängen. Ein Unterricht, in dem konsequent \emph{Flächen} (Dreiecksflächen, Rechtecksflächen, ...) zerlegt und wieder zusammengesetzt werden und nicht \emph{Objekte} (Dreiecke, Rechtecke, ...), fördert nicht die Vorstellung, Flächen seien unteilbare Ganzheiten. \smallskip

Fassen wir zusammen. Ein das Theorie -- Element modellierendes Arbeitsmittel setzt sich zwei Ziele:\\
\noindent --- Flächen lassen sich auf unterschiedliche Weise zerlegen.\\
\noindent --- Setzt man aus den Bausteinen, in die eine Fläche zerlegt wird, eine Fläche anderer Form zusammen, so hat diese den gleichen Flächeninhalt, da dieselben Bausteine verwendet und nur anders angeordnet werden.\smallskip

Damit soll den empirisch erwiesenen Defiziten, daß Maßquadrate nicht teilbar sind, und zerlegungsgleiche Figuren nicht gleichen Flächeninhalt haben, begegnet werden. Ob man diese Defizite damit ausräumt, kann nur die Empirie beantworten. \smallskip

Selbst wenn die genannten Ziele des Arbeitsmittels erreicht werden, ist damit die premeasurement -- Phase nicht abgeschlossen. Lynne Outhred und Michael Mitchelmore betonen die Wichtigkeit, die Zeichnungen für die Bestimmung des Flächeninhaltes haben. Die zum Verständnis der \emph{Flächenformel} erforderliche Einsicht, daß man einem Rechteck eine Zeilen -- Spalten -- Struktur aufprägen kann, ist erst verinnerlicht, wenn diese Struktur zeichnerisch beherrscht wird. Wird sie nur mit Bausteinen gelegt, bleibt sie für den Schüler an das Material gebunden. Wie die Untersuchung von Outhred und Mitchelmore zeigt, ist im Grundschulalter jeder zweite Schüler außerstande, ein mit Bausteinen gelegtes Muster zu zeichnen, sogar es abzuzeichnen [Outhred und Mitchelmore 1992].

Sieht man die Einsichten, die ein das Theorie -- Element T$_{FV}$ modellierendes Arbeitsmittel auf enaktiver Ebene vermitteln kann, als erste Stufe des Lernprozesses an, der zum Verständnis der Flächeninhaltsbestimmung führt, so wären folglich auf  einer zweiten Stufe die Aktivitäten, die das Arbeitsmittel erlaubt, auf ikonischer (zeichnerischer) Ebene zu wiederholen. Beherrschen die Schüler die durch das Theorie -- Element T$_{FV}$ vermittelten Einsichten auf der enaktiven Ebene und sind sie in der Lage, diese auf die ikonische Eben zu übertragen, so ist die premeasurement -- Phase abgeschlossen.

Da die zeichnerischen Fähigkeiten der Schüler und das Verständnis der Homogenität einer Fläche für das Verständnis der Flächeninhaltsbestimmung große Bedeutung haben, ist dem Unterricht damit eine zeitliche Vorgabe gegeben, die es einzuhalten gilt, will man das Thema erfolgreich angehen. \bigskip

Bevor wir uns einer \emph{Empirischen Theorie der Flächeninhaltsbestimmung} (der betrachteten Polygonflächen) zuwenden, also einer Theorie zur numerischen Bestimmung von Flächeninhalten, halten wir noch einmal fest, daß die premeasurement -- Phase durch ein Theorie -- Element beschrieben wurde, das ausschließlich empirische Objekte zum Gegenstand hat. Denn auch Zeichenblattfiguren sind für den Schüler reale Objekte. Man kann sie ja ausschneiden und wie Plättchen handhaben. \smallskip

Als erstes formulieren wir eine wesentliche Voraussetzung, auf die die folgende Theorie zurückgreifen wird. \\
\noindent Vorausgesetzt wird, daß die Schüler über einen angemessenen Maßzahlbegriff verfügen, daß sie ein relationales und nicht nur ein instrumentelles Verständnis des Messens von Längen haben [Outhred und Mitchelmore 2000]. Dazu  sind drei Einsichten von besonderer Bedeutung: \\
\noindent --- Um die Länge eines Objektes (Tischkante, gespannter Faden o.ä.) zu ermitteln wird dieses mit einem frei wählbaren (kleineren) Objekt lückenlos und überschneidungsfrei \emph{ausgelegt}. \\
\noindent --- Das numerische Ergebnis der Messung ist abhängig von der Größe des gewählten Vergleichsobjektes. \\
\noindent --- Die Wahl eines Realisanten einer einheitlichen Maßeinheit$(10^{n}$-ter Teil (n $\in \mathbb{N}_{0})$ des Urmeters) dient ausschließlich dazu, die Vergleichbarkeit von Messungen zu vereinfachen. \\

Daß diese Einsichten keineswegs mit einem nur instrumentellen Verständnis der Längenmessung einhergehen, belegen die Ausführungen in [Kamii und Kysh 2000]. \smallskip 

Wir betonen nochmals, daß dem Schüler bewußt sein muß, daß beim linearen Messen \emph{ausgelegt} und nicht nur gezählt wird. Deshalb haben wir bei der Besprechung des Arbeitsmittels dem Auslegen von Gebilden mit Bausteinen solche Bedeutung beigemessen. Es soll die Schüler an das Auslegen beim linearen Messen erinnern. \smallskip 

Neben der genannten Voraussetzung spielen bei der \emph{Empirischen Theorie der Flächeninhaltsbestimmung} Vereinbarungen eine wesentliche Rolle, deren Gründe für den Schüler nicht unmittelbar ersichtlich sein dürften.

Wer über das Theorie -- Element T$_{FV}$ verfügt, kann Flächen der Größe nach vergleichen, ohne Zahlen heranzuziehen (vgl. M(T$_{FV}$) (vii)), wobei wir davon ausgehen, daß beim Nachweis von U~zgl~V dieselben Bausteine benutzt werden, um die Einsicht zu fördern, daß die Zerlegungsgleichheit die Flächengleichheit von U und V impliziert. 

Will man zahlenmäßige Vergleiche vornehmen, so empfiehlt es sich, beide Flächen mit kongruenten Bausteinen auszulegen. Wie Bernard Héraud zeigte, neigen Schüler überwiegend dazu, zum Auslegen kongruente Bausteine zu wählen, dabei bevorzugt solche, die der auszulegenden Fläche ähnlich sind, d.h. eine rechteckige Fläche wird mit rechteckigen Bausteinen ausgelegt und eine dreieckige mit dreieckigen [Héraud 1987].  

Im Theorie -- Element T$_{FV}$ haben wir die Zerlegungsgleichheit eingeführt, die es erlaubt, alle in Frage kommenden Gebilde in inhaltsgleiche rechteckige Gebilde zu überführen (vgl. M(T$_{FV}$) (viii)). Es bleibt aber die Frage, weshalb man die zahlenmäßige Inhaltsangabe auf Quadraten aufbaut. Diese konventionelle Festsetzung kann nicht aus einem rein geometrischen Argument abgeleitet werden. Unserer Auffassung nach ist der einzige Weg, dem Schüler gegenüber intellektuell ehrlich vorzugehen, die anvisierten Rechtecksflächen bewußt in die Überlegung einzubeziehen.

Zu dem vorausgesetzten relationalen Verständnis des linearen Messens  gehört die Einsicht, daß man zwar die Einheit der Längenmessung frei wählen kann, dann aber aus Gründen der \emph{Zweckmäßigkeit} und nicht aus \emph{logischen} Gründen alle linearen Messungen mit dieser Einheit durchführt. Ferner war ein Ziel des das Theorie -- Element modellierenden Arbeitsmittels die Einsicht, daß Flächen verschiedener Form gleiche Größe haben können, und daß man dies durch Auslegen überprüfen kann.

Hinzukommt, wie Konstantinos Zacharos beim Vergleich der Inhaltsbestimmung durch Auslegen bzw. die \emph{Flächenformel} --- die ja auf lineares Messen zurückgreift --- feststellte, daß Schüler beim Messen erfolgreicher sind, wenn die Einheit des "Meßwerkzeugs" die gleiche Dimension hat wie das auszumessende Gebilde, was ebenfalls für das Auslegen spricht [Zacharos 2006].

Da die Zerlegungsgleichheit erlaubt, jede Polygonfläche in eine rechteckige Fläche zu überführen, reduziert sich das Problem --- unter Berücksichtigung des oben zitierten Ergebnisses von Héraud --- auf die Wahl einer geeigneten rechteckigen Fläche.

In [Battista u. a. 1998] und [Outhred und Mitchelmore 2000] wird übereinstimmend betont, daß die Schüler die Fläche eines Rechtecks in geeigneter Weise strukturieren müssen, bevor sie in der Lage sind, ihren Inhalt zu bestimmen. Die Struktur muß dazu verinnerlicht werden, was nach [Outhred und Mitchelmore 1992] bedeutet, daß sie auch zeichnerisch beherrscht wird (s.o.).

Damit ist aber nach wie vor offen, weshalb ein Schüler als Strukturierung ein aus Quadraten aufgebautes Netz wählen soll. \\

Wie schon angesprochen, ist das Messen durch zwei Spezifika charakterisierbar: \\
\noindent --- die Wahl der Einheit, \\
\noindent --- das Auslegen des zu messenden Objektes mit dieser Einheit.

Was das Auslegen eine Fläche und damit die \emph{additive} Bestimmung des Flächeninhaltes betrifft, so wurde diese Komponente des Messens in die \emph{Empirische Theorie des Flächenvergleichs} verlagert, und wir sind der Auffassung, daß sie dort mit der erforderlichen Verständnistiefe behandelt werden kann. Was verbleibt, ist die Wahl der Einheit.

\noindent Naheliegend für Schüler sind 

\noindent --- eine zweidimensionale Einheit zu wählen [Zacharos 2006],

\noindent --- in gleicher Weise wie beim linearen Messen für jede Messung \emph{genau eine} Einheit zu wählen [Héraud 1987].\\

Man darf davon ausgehen, daß der Schüler bislang nur das lineare Messen in sonderlicher Breite und Tiefe kennengelernt hat. Was er noch nicht kennen dürfte, ist, daß es zweckmäßig sein kann, Maßeinheiten voneinander abhängig zu machen. Die Wahl eines Einheitsquadrates als Maßeinheit für die Flächenmessung ist damit nicht \emph{logisch} bedingt, sondern folgt einem Zweckmäßigkeitsgesichtspunkt. Nur wenn dieses im Unterricht hinreichend betont wird, vermag der Schüler einen Zusammenhang zwischen der Flächenformel und dem Auslegen mit Einheitsquadraten herzustellen, dessen Fehlen beklagt wird [Simon und Blume 1994], [Outhred und Mitchelmore 2000].

Vom Schüler zu erwarten, daß er ein Einheitsquadrat als Einheit wählt, ohne den leitenden Gesichtspunkt zu kennen, erscheint uns etwas weit hergeholt. Nennt man ihm aber diesen,  so liegt die Wahl des Einheitsquadrates nahe, und eine Rechteckfläche dann mit Quadraten auszulegen ist eine fast zwangsläufige Folgerung. Das Quadratnetz als Zeilen -- Spalten -- Struktur zu sehen ist für Kinder zwar keineswegs naheliegend [Battista u. a. 1998], dürfte aber erleichtert werden, wenn man sie darauf hinweist, daß das Maß der Fläche mit Hilfe der Längen der Rechteckseiten bestimmt werden soll. Es gibt allerdings auch Untersuchungen, die zeigen, daß diese Hilfe nicht generell notwendig ist. Um den Flächeninhalt des Rechtecks durch Auslegen zu bestimmen gliederten die Kinder die erforderlichen Plättchen in Zeilen / Spalten und ermittelten die Gesamtzahl durch iterative Addition oder Multiplikation [Nunes u. a. 1993], [Bonotto 2003]. \bigskip

Wir kommen nun zur Formulierung eines Theorie -- Elementes das die \emph{Empirischen Theorie der Flächeninhaltsbestimmung} formal beschreibt. Wir bezeichnen es mit T$_{FI}$ (FI für "Flächeninhalt"). Es handelt von Polygonflächen --- realisiert durch die Gebilde der \emph{Empirischen Theorie des Flächenvergleichs} sowohl zusammengefügt aus Bausteinen wie dargestellt in Zeichnungen --- und der Zuweisung von Zahlenwerten zu diesen Objekten. Der Charakter dieser Objekte stellt sicher, daß es sich um eine empirische Theorie handelt. \medskip 

Wir beginnen mit der Formulierung der \emph{didaktischen Postulate}, denen die Theorie zu genügen hat. \smallskip \\ 
\noindent (P 1)  Der Zweck, zu dem ein Schüler die Bestimmung von Flächeninhalten erlernt, ist die Bewältigung von Alltagsproblemen. \medskip \\
\noindent(P 2)  Es werden Maßzahlen für die Inhalte von Polygonflächen eingeführt, und es wird gezeigt, wie man sie bestimmt. \smallskip \\
\noindent (P 3)  Systematische Voraussetzung ist ein relationales Verständnis des linearen Messens. \medskip 

\noindent Nun die Angabe des Theorie -- Elementes T$_{FI}$. 

Vorausgesetzt wird die Beherrschung eines Theorie -- Elementes, das angibt, was es heißt, über einen Maßzahlbegriff für lineares Messen zu verfügen, der auf der freien Wahl einer Einheit und dem Auslegen der zu messenden Objekte beruht.  \smallskip

Die \emph{partiellen Modelle} bestehen je aus einem Modell des Theorie -- Elementes T$_{FV}$, ergänzt um ein Modell eines vorausgesetzten Theorie -- Elementes des linearen Messens. \smallskip

Die \emph{potentiellen Modelle} enthalten als theoretischen Begriff die (spätere) Flächeninhaltsfunktion $\mu$. Da wegen der in T$_{FV}$ geforderten Zerlegungsgleichheit jede Polygonfläche in ein inhaltsgleiches Rechteck verwandelt werden kann (vgl. M(T$_{FV})$ (viii)), können wir $\mu$ auf Rechtecksflächen beschränken --- auf dieser ersten Stufe auf solche mit ganzzahligen Seitenlängen. Hat ein Rechteck die Seitenlängen me und ne (e für die Längeneinheit), so bezeichnen wir es mit ($\overline{me,ne}$). Schreiben wir $\mathbb{N}_{e}$ für $\{me \lvert m \in \mathbb{N}\}$ (und entsprechend $\mathbb{N}_{e^2}$), so läßt sich $\mu$ auffassen als \\
\begin{center}
$\mu \subseteq \big(\overline {\mathbb{N}_{e} \times \mathbb{N}_{e}},\mathbb{N}_{e^2} \big)$ \\ \bigskip
\end{center}
Die \emph{Modelle} enthalten schließlich die (übliche) Definition der Flächeninhaltsfunktion.
\begin{center}
\begin{tabular}{l l}
 $\mu $:& $ \overline{\mathbb{N}_{e} \times \mathbb{N}_{e}} \quad \rightarrow \quad \mathbb{N}_{e^2}$ \\ 
& $\big(\overline{me,ne}\big)\;  \mapsto \;  \big(m \cdot n \big)e^2$\medskip 
\end{tabular}
\end{center}

Was die Adäquatheit dieser Definition betrifft, so haben wir schon darauf verwiesen, daß wir es für notwendig halten, die Zuordnung $\big(\overline{1e,1e}\big) \: \mapsto \: 1e^2$ vorzugeben. Damit wird verständlich, daß $1e^2$ eine Einheit zum Messen von Flächen ist, und wie sich $1e^2$ aus dem Längenmaß 1e ableitet. Die Gefahr, daß Schüler keine Verbindung herstellen zwischen einer Längeneinheit und der ihr entsprechenden Flächeneinheit [Simon und Blume 1994], [Bonotto 2003] dürfte man damit begegnen können. \smallskip

Die \emph{Querverbindung} schließlich stellt sicher, daß nur solche Mengen potentieller Modelle in das Theorie -- Element $T_{FI}$ aufgenommen werden, bei denen $\mu$ für jedes Element bei zerlegungsgleichen Flächen gleiche Werte liefert. \bigskip

Was unterscheidet nun das hier vorgestellte Konzept von der in den meisten Schulbüchern gewählten Vorgehensweise? Wir halten zwei Punkte für entscheidend:

\noindent --- Wir ziehen eine breit aufgefächerte \emph{Empirische Theorie des Flächenvergleichs} vor, um alle für die Flächeninhaltsbestimmung wichtigen Einsichten auf enaktiver und \emph{anschließend} auf ikonischer (zeichnerischer) Ebene vorzubereiten, wie es die Empiriker fordern (vgl. (Osborne o.J. [1976])).

\noindent --- Da in der \emph{Empirischen Theorie der Flächeninhaltsbestimmung} eine Entscheidung eine wichtige Rolle spielt, die nach dem Gesichtspunkt der Zweckmäßigkeit getroffen wird --- die Wahl eines Quadrates als Maßeinheit ---, erhält sie einen normativen Aspekt, wenn man diese Entscheidung als Sollbestimmung formuliert. Darin kommt deutlich zum Ausdruck, daß das hier vorgelegte Konzept von den schulüblichen Vorgehensweisen abweicht. \bigskip

Welche Konsequenzen hat das hier vorgelegte Konzept?

\noindent --- Der zeitliche Umfang der Behandlung des Themas dürfte deutlich zunehmen. 

\noindent --- Schon die Behandlung der \emph{Empirischen Theorie des Flächenvergleichs} auf der ikonischen Ebene muß auf einer weiterführenden Schule erfolgen, da die kognitiven Leistungen und zeichnerischen Fähigkeiten, die der Schüler zu erbringen hat, im Grundschulalter nicht erbracht werden können. \medskip

Man sieht, daß die Konsequenzen, die die empirischen Ergebnisse einfordern, die Auswahl der Unterrichtsinhalte an einem neuralgischen Punkt treffen, dem Zeitfaktor. Aber die vorgelegte Konzeption bietet die Möglichkeit, begründet abzuwägen, in welchem Umfang man den Ergebnissen der Empiriker Rechnung tragen will. Lediglich die vorhandenen Defizite zu konstatieren führt nicht weiter. \bigskip 

Wie man dieses zweite Beispiel bewertet, hängt davon ab, in welchen Rahmen man sich empirische Ergebnisse eingearbeitet wünscht. Bevorzugt man sie nur lokal organisiert, so sind Konzepte wie das hier vorgestellte überflüssig. Wünscht man sie dagegen stärker zu strukturieren, so stellt die strukturalistische Darstellungsform der empirischen Theorie ein zweckmäßiges Rüstzeug bereit. Dann wird man auch auf die formale Darstellung nicht verzichten. Es wäre absurd, sie durch eine umgangssprachliche Beschreibung zu ersetzen, die wesentlich umfangreicher, unpräziser und weniger detailliert wäre. Die formale Darstellung weist zudem die Punkte aus, in denen die empirische Rechtfertigung der Theorie noch Lücken hat --- darauf wurde schon verwiesen.
\bigskip
 
Wir haben einleitend darauf hingewiesen, daß eine empirische Theorie deskriptiv oder präskriptiv / normativ gelesen werden kann. Die hier vorgelegten Theorie - Elemente des Flächenvergleichs und der Flächeninhaltsbestimmung, die eine Vielzahl empirischer Ergebnisse aufgreifen, um als Basis einer Unterrichtskonzeption dienen zu können, sind folglich präskriptiv / normativ zu lesen.
\vspace{0,5 cm} 

Unter Verwendung des strukturalistischen Begriffsrahmens haben wir an zwei äußerst unterschiedlichen Beispielen gezeigt, wie sich das Begriffssystem nutzen läßt, mathematikdidaktisch relevante Ergebnisse zu systematisieren. Man kann natürlich fragen, weshalb zur Darstellung der beiden Beispiele der strukturalistische Begriffsrahmen gewählt wurde. Dazu wäre zu fragen, weshalb es nicht legitim sein soll, einen vorliegenden begrifflichen Rahmen zu verwenden. Seine Schlüssigkeit und praktische Anwendbarkeit kann der Leser --- zumindest im Ansatz ---  an Hand des vorgelegten Textes beurteilen. Dazu dient dieser. Den Verfassern ist jedenfalls keine geeignete alternative Darstellungsform bekannt. Wenn Kinder im Mathematikunterricht empirische Theorien erlernen --- so die These der Autoren, die nicht nur durch mathematikdidaktische Untersuchungen gestützt wird sondern auch durch kognitionspsychologische (vgl. [Gopnik und Meltzoff 1997]) ---, dann ist es ein angemessenes wissenschaftliches Vorgehen, einschlägige Darstellungsmittel aus anderen wissenschaftlichen Disziplinen --- hier der Wissenschaftstheorie -- zu benutzen und auf ihre Tragfähigkeit zu testen. Dies geschieht auf unterschiedliche Weise in den beiden Teilen des vorliegenden Beitrages.   
\vspace{2cm}

\selectlanguage{german}
\nocite{*}
\printbibliography

\end{document}